\documentclass[12pt]{article}
  \usepackage{amsfonts}
  \usepackage[square,sort&compress,comma,numbers]{natbib} 
   \def\R{\mathbb{R}}

   \def\1{{\rm I\mskip -10.5mu 1}} 
   
   \def\e{{\varepsilon}}

   \def\cC{{\cal C}}
   \def\cD{{\cal D}}

   \def\dist{\mathop{\rm dist}\nolimits}
   
   \def\div{\mathop{\rm div}\nolimits}

   \def\no{\noindent}
   \def\proof{\mbox {{\underline {\sf Proof}} \hspace{2mm}}}
   \def\qed{{\hfill {\em q.e.d.}\\\vspace{1mm}}}
   
   \newcommand{\beq}{\begin{equation}}
   \newcommand{\eeq}{\end{equation}}

\newtheorem{df}{Definition}[section]
\newtheorem{prop}[df]{Proposition}
\newtheorem{lemma}[df]{Lemma}
\newtheorem{teo}[df]{Theorem}
\newtheorem{rem}[df]{Remark}

\newtheorem{cor}[df]{Corollary}

 \newcommand{\sezione}[1]{\section{#1}\setcounter{equation}{0}}

\begin{document}


   \title{Uniqueness of solutions for nonlinear Dirichlet problems
     with supercritical growth
}


  \maketitle


 \vspace{5mm}

\begin{center}

{ {\bf Riccardo MOLLE$^a$,\quad Donato PASSASEO$^b$}}

\vspace{5mm}

{\em
${\phantom{1}}^a$Dipartimento di Matematica,
Universit\`a di Roma ``Tor Vergata'',\linebreak
Via della Ricerca Scientifica n. 1,
00133 Roma, Italy.}

\vspace{2mm}

{\em
${\phantom{1}}^b$Dipartimento di Matematica ``E. De Giorgi'',
  Universit\`a di Lecce,\linebreak 
P.O. Box 193, 73100 Lecce, Italy.
}
\end{center}

\vspace{5mm}


{\small {\sc \noindent \ \ Abstract.} - 
\footnote{ {\em E-mail address:} molle@mat.uniroma2.it (R. Molle).}
We are concerned with Dirichlet problems of the form
$$
\div(|D u|^{p-2}Du)+f(u)=0\ \mbox{ in }\Omega,\qquad 
u=0\ \mbox{ on }\partial\Omega,
$$
where $\Omega$ is a bounded domain of $\R^n$, $n\ge 2$, $1<p<n$ and
$f$ is a continuous function with supercritical growth from the
viewpoint of the Sobolev embedding.

\no In particular, if $n=2$ and $\gamma:[a,b]\to\R^2$ is a smooth
curve such that $\gamma(t_1)\neq\gamma(t_2)$ for $t_1\neq t_2$, we
prove that, for $\e>0$ small enough, there exists a unique solution of
the Dirichlet problem in the domain
$\Omega=\Omega^\Gamma_\e=\{(x_1,x_2)\in\R^2\ :\
\dist\big((x_1,x_2),\Gamma\big)<\e\}$, where $\Gamma=\{\gamma(t)\ :\
t\in[a,b]\}$. 

\no Moreover, we extend this uniqueness result to the case where $n>2$
and $\Omega$ is, for example, a domain of the type
$$
\Omega=\widetilde\Omega^\Gamma_{\e,s}=\{(x_1,x_2,y)\ :\
(x_1,x_2)\in\Omega^\Gamma_\e, \ y\in\R^{n-2},\ |y|<s\}.
$$

\vspace{3mm}


{\em  \noindent \ \ MSC:}  35J20; 35J60; 35J65.

\vspace{1mm}

{\em  \noindent \ \  Keywords:} 
   Supercritical Dirichlet problems, contractible domains,
   nonexistence of solutions.
}


\sezione{Introduction}


In this paper we deal with nonlinear Dirichlet problems of the form 
\beq
\label{P}
\div(|D u|^{p-2}Du)+f(u)=0\ \mbox{ in }\Omega,\qquad u=0\ \mbox{ on
}\partial \Omega, 
\eeq
where $\Omega$ is a bounded domain of $\R^n$, $n\ge 2$, $1<p<n$ and
$f$ is a continuous function that, for a suitable $q>{np\over n-p}$,
satisfies the condition
\beq
\label{f}
t\, f(t)\ge q\int_0^tf(\tau)d\tau\ge 0\qquad\forall t\in\R
\eeq
(this means that $f$ has a supercritical growth from the viewpoint
of the Sobolev embedding $H^{1,p}_0(\Omega) \hookrightarrow
L^q(\Omega)$).

\no It is well known that the existence of
nontrivial solutions for problem (\ref{P}) is strictly related to the
shape of $\Omega$ (see \cite{B,BN}).
For example, if $\Omega$ is an annulus there exist infinitely many
solutions (see f.i. \cite{KW}), while if $\Omega$ is star-shaped the
problem has only the trivial solution $u\equiv 0$ as a consequence of a Pohozaev
type identity (see \cite{Po}).

\no In this paper our aim is to show that this uniqueness result may
be extended to some bounded contractible non star-shaped domains $\Omega$
that can be very different from the star-shaped ones and even arbitrarily
close to non contractible domains.

\no If $n=2$, we construct these domains in the following way.
Given a smooth curve $\gamma:[a,b]\to\R^2$ such that $\gamma'(t)\neq
0$ $\forall t\in[a,b]$ and $\gamma(t_1)\neq\gamma(t_2)$ for $t_1\neq
t_2$, we set $\Gamma=\{\gamma(t)\ :\ t\in[a,b]\}$ and, for all $\e>0$,
we consider the domain $\Omega=\Omega^\Gamma_\e$ defined by
\beq
\label{Oe}
\Omega^\Gamma_\e=\{(x_1,x_2)\in\R^2\ :\ \dist\big((x_1,x_2),\Gamma\big)<\e\}.
\eeq
We prove that, for $\e>0$ small enough, the Dirichlet problem (\ref{P})
with $\Omega=\Omega^\Gamma_\e$ has only the trivial solution $u\equiv
0$ (see Theorem \ref{T2.0}).

\no It is clear that $\Omega^\Gamma_\e$ is contractible for $\e>0$ small enough.
Moreover, it is not star-shaped (unless $\Gamma$ is a segment of a
stright line) and it may be arbitrarily close to a non contractible
domain (because $|\gamma(a)-\gamma(b)|$ may be arbitrarily small).
This fact (as we pointed out also in \cite{MPmed}) seems to suggest
that for $n=2$ we have existence of nontrivial solutions when $\Omega$
is not contractible and nonexistence when $\Omega$ is contractible. 

\no For $n>2$ the situation is more complex because there exist
contractible domains $\Omega$, even arbitrarily close to star-shaped
domains, such that the problem has nontrivial solutions (see for example 
\cite{BaCo,B,BN,CCL,Co,D88,DZ,Di,KW,MPcras02,MPcras2002,MPaihp04,MPcvpde06,MPaihp06,MPT,Pmm89,Pl92,P4092,Pjfa93,P93,P94,Pdie95,Pna95,Ptmna96,Pd98,Po,PT,R},
where the effect of the domain shape on the number of solutions is
studied, answering some well-known questions posed by Brezis,
Nirenberg, Rabinowitz, etc.).

\no However, also for $n>2$ we can obtain uniqueness results in
bounded, contractible, non star-shaped domains $\Omega$ of $\R^n$,
arbitrarily close to non contractible domains.
For example, we can consider domains of the type
\beq
\label{Otilde}
\widetilde\Omega^\Gamma_{\e,s}=\{(x_1,x_2,y)\ :\
(x_1,x_2)\in\Omega^\Gamma_\e, \ y\in\R^{n-2},\ |y|<s\}
\eeq
and prove that, for $\e>0$ small enough and $s>0$, the Dirichlet
problem (\ref{P}) with $\Omega=\widetilde\Omega^\Gamma_{\e,s}$ has
only the trivial solution $u\equiv 0$ (see Theorem \ref{T3.0}).


\sezione{Uniqueness result in the case $n=2$}


The main result in the case $n=2$ is presented in the following
theorem. 

\begin{teo}
\label{T2.0}
Assume that $\gamma\in\cC^3([a,b],\R^2)$, $\gamma'(t)\neq 0$ $\forall
t\in[a,b]$ and $\gamma(t_1)\neq\gamma(t_2)$ for $t_1\neq t_2$.
Let $\Omega=\Omega^\Gamma_\e$ be defined as in (\ref{Oe}).
Moreover, assume that $1<p<2$ and there exists $q>{2p\over 2-p}$ such
that condition (\ref{f}) holds.
Then, there exists $\bar\e>0$ such that the Dirichlet problem
(\ref{P}) has only the solution $u\equiv 0$ for all $\e\in(0,\bar\e)$.
\end{teo}

\no The proof requires some preliminary results.

\no Since $\gamma$ is a smooth curve such that $\gamma'(t)\neq 0$
$\forall t\in[a,b]$ and $\gamma(t_1)\neq\gamma(t_2)$ for $t_1\neq t_2$,
there exists $\bar\e_1>0$ such that 
$\Omega^\Gamma_{\bar\e_1}$ is a contractible domain and, for all
$(x_1,x_2)\in\overline\Omega^\Gamma_{\bar\e_1}$, there exists a unique $t\in[a,b]$
satisfying $\dist\big((x_1,x_2),\gamma(t)\big)=\dist \big(
(x_1,x_2),\Gamma\big)$.
Without any loss of generality, we can assume in addition that
$a\le0\le b$ and $|\gamma'(t)|=1$ $\forall t\in[a,b]$.
Let us denote by
$\gamma_{\bar\e_1}:[a-{\bar\e_1},b+{\bar\e_1}]\to\R^2$ the curve such
that  
\beq
\gamma_{\bar\e_1}(t)=\gamma(t)\quad\forall t\in[a,b],\qquad
\gamma'_{\bar\e_1}(t)=\gamma'(a)\quad\forall t\le a,\qquad 
\gamma'_{\bar\e_1}(t)=\gamma'(b)\quad\forall t\ge b.
\eeq
Moreover, let us set
\beq
T(t)=\gamma'_{\bar\e_1}(t)=\left( 
\gamma'_{{\bar\e_1},1}(t), \gamma'_{{\bar\e_1},2}(t) \right),
\quad
N(t)=\left(-\gamma'_{{\bar\e_1},2}(t), \gamma'_{{\bar\e_1},1}(t) \right)\qquad \forall t\in[a-\bar\e_1,b+\bar\e_1].
\eeq
Then, for all $(x_1,x_2)$ in $\overline\Omega^\Gamma_{\bar\e_1}$  there exists a unique pair
$(t,r)\in\R^2$ such that 
\beq
t\in[a-{\bar\e_1}, b+{\bar\e_1}],\quad r\in[-{\bar\e_1} ,
{\bar\e_1}]\quad\mbox{ and }\quad
(x_1,x_2)=\gamma_{\bar\e_1}(t)+rN(t).
\eeq
Since $\Gamma\in\cC^3([a,b],\R^2)$, we can consider in
$\overline\Omega^\Gamma_{\bar\e_1}$ the vector field
$v=(v_1,v_2)\in\cC^1(\overline\Omega^\Gamma_{\bar\e_1},\R^2)$ defined
by
\beq
\label{vf}
v(\gamma_{\bar\e_1}(t)+rN(t))=tT(t)[1-r\,
\gamma_{\bar\e_1}''(t)\cdot N(t)]+rN(t)\quad\forall t\in[a-{\bar\e_1},
b+{\bar\e_1}],\ \forall r\in[-{\bar\e_1},{\bar\e_1} ].
\eeq
In next lemma we describe the main properties of the vector field
$v$.
\begin{lemma}
\label{L2.1}
If the curve $\gamma$ satisfies all the above required assumptions,
the vector field $v\in\cC^1(\overline\Omega^\Gamma_{\bar\e_1},\R^2 )$ defined in
(\ref{vf}) satisfies
\begin{itemize}
\item[{\em a)}]$v\cdot\nu>0$ on $\partial\Omega^\Gamma_\e$
  $\forall\e\in(0,\bar\e_1)$; 
\item[{\em b)}] $\div
  v(\gamma_{\bar\e_1}(t)+rN(t))=2-{r[t\gamma_{\bar\e_1}''(t)\cdot
    N(t)]'\over 1-r\gamma_{\bar\e_1}''(t)\cdot
    N(t)}$ \quad $\forall t\in[a-{\bar\e_1},
b+{\bar\e_1}],\ \forall r\in[-{\bar\e_1},{\bar\e_1} ]$,
\item[{\em c)}]
$dv(\gamma_{\bar\e_1}(t)+rN(t))[\xi]\cdot\xi=
\left[1-{r[t\gamma_{\bar\e_1}''(t)\cdot
    N(t)]'\over 1-r\gamma_{\bar\e_1}''(t)\cdot
    N(t)}\right]\,\xi^2_T(t)+\xi^2_N(t)$ \quad $\forall\xi\in\R^2$, $\forall
t\in[a-{\bar\e_1}, 
b+{\bar\e_1}],\ \forall r\in[-{\bar\e_1},{\bar\e_1} ]$,
where 
$\xi_T(t)=\xi\cdot T(t)$
and 
$\xi_N(t)=\xi\cdot N(t)$
\end{itemize}
(here $\nu$ denotes the outward normal to
$\partial\overline\Omega^\Gamma_{\e}$ and
$dv[\xi]=\sum\limits_{i=1}^2\xi_i\, D_iv$ $\forall\xi=(\xi_1,\xi_2)\in\R^2$).
\end{lemma}

\proof
Taking into account that $a\le 0\le b$, as we have assumed, property
$(a)$ is a direct consequence of the choice of $\bar\e_1$ and the definition of $\Omega^\Gamma_\e$
and $v$.

\no In order to prove $(b)$ and $(c)$, notice that for all
$t\in[a-\bar\e_1,b+\bar\e_1]$ and $r\in[-\bar\e_1,\bar\e_1]$ we have
\beq
\label{N}
dv(\gamma_{\bar\e_1}(t)+rN(t))[N(t)] 
= 
-t\gamma_{\bar\e_1}''(t)\cdot N(t)\, T(t)+N(t),
\eeq
as one can verify by direct computation, and
\begin{eqnarray}
\label{T}
dv(\gamma_{\bar\e_1}(t)+rN(t))[T(t)] 
& = &
[1-r\gamma_{\bar\e_1}''(t)\cdot N(t)]^{-1}{\partial \over\partial t}
v(\gamma_{\bar\e_1}(t)+rN(t))
\\
\nonumber
& = & \left[1-{r[t\gamma_{\bar\e_1}''(t)\cdot
    N(t)]'\over 1-r\gamma_{\bar\e_1}''(t)\cdot
    N(t)}\right] \, T(t)+t\gamma_{\bar\e_1}''(t)\cdot N(t)\, N(t)
\end{eqnarray}
because
\beq
{\partial\over\partial
  t}[\gamma_{\bar\e_1}(t)+rN(t)]=[1-r\gamma_{\bar\e_1}''(t)\cdot
N(t)]T(t)
\eeq
(notice that $[1-r\gamma_{\bar\e_1}''(t)\cdot
N(t)]>0$ $\forall t\in[a-\bar\e_1,b+\bar\e_1]$ if $\bar\e_1$ is small enough).

\no Then $(b)$ and $(c)$ follow from (\ref{T}) and (\ref{N}).

\qed

\begin{lemma}
\label{L2.2}
If the curve $\gamma$ satisfies all the above required assumptions, we have
\beq
\label{lim}
\lim_{\e\to 0}\max\left\{
\left| {r[t\gamma_{\bar\e_1}''(t)\cdot
    N(t)]'\over 1-r\gamma_{\bar\e_1}''(t)\cdot
    N(t)}\right|\ :\ -\e\le r\le \e,\ t\in [a-\e,b+\e]
\right\}=0.
\eeq
\end{lemma}

\proof Since $\gamma\in\cC^3([a,b],\R^2)$, the maximum in (\ref{lim})
is achieved for all $\e\in(0,\bar\e_1)$.
If it is achieved on the pair $(t_\e,r_\e)$, we have
$\lim\limits_{\e\to 0}t_\e=\bar t$ (up to a subsequence) for a
suitable $\bar t\in[a,b]$ and $\lim\limits_{\e\to 0}r_\e=0$ (because
$|r_\e|\le\e$). 
Then (\ref{lim}) follows easily.

\qed

\no The following lemma generalizes Pohozaev identity.

\begin{lemma}
\label{L2.3}
Assume that, for all $\e\in(0,\bar\e_1)$, $u_\e$ is a solution of the
Dirichlet problem
\beq
\label{eq}
\div(|Du_\e|^{p-2}Du_\e)+f(u_\e)=0\quad\mbox{ in }\Omega^\Gamma_\e,
\qquad u_\e=0\ \mbox{ on } \partial \Omega^\Gamma_\e
\eeq
and consider a vector field
$v=(v_1,v_2)\in\cC^1(\overline\Omega_\e^\Gamma,\R^2)$.

\no Then the following integral identity holds:
\beq
\label{id}
\left(1-{1\over
    p}\right)\int_{\partial\Omega^\Gamma_\e}|Du_\e|^pv\cdot\nu\,
d\sigma =\phantom{***********************}
\eeq
$$
\phantom{*********}
=\int_{\Omega^\Gamma_\e}|Du_\e|^{p-2}dv[Du_\e]\cdot Du_\e dx+
\int_{\Omega^\Gamma_\e}\div v\left[F(u_\e)-{1\over p}|Du_\e|^p\right]dx,
$$
where $F(t)=\int_0^t f(\tau)d\tau$ $\forall t\in\R$.
\end{lemma}

\proof
From (\ref{eq}) we infer that
\beq
\int_{\Omega^\Gamma_\e}\div(|Du_\e|^{p-2}Du_\e)\,v \cdot Du_\e dx+
\int_{\Omega^\Gamma_\e}f(u_\e)\, v\cdot Du_\e dx=0
\eeq
which implies
\beq
\label{e1}
\int_{\partial\Omega^\Gamma_\e}
|Du_\e|^{p-2} (Du_\e\cdot\nu)(v\cdot Du_\e) d\sigma =
\int_{\Omega^\Gamma_\e}
|Du_\e|^{p-2} Du_\e\cdot D(v\cdot Du_\e) dx 
-\int_{\Omega^\Gamma_\e}f(u_\e)\, v\cdot Du_\e dx.
\eeq
Since $u\equiv 0$ on $\partial\Omega^\Gamma_\e$, we have $D
u_\e=(Du_\e\cdot\nu)\, \nu$ and, taking into account the definition of
$F$, $F(u)=0$ on $\partial\Omega^\Gamma_\e$.
As a consequence, we obtain
\beq
\label{e2}
\int_{\partial\Omega^\Gamma_\e}
|Du_\e|^{p-2} (Du_\e\cdot\nu)(v\cdot Du_\e) d\sigma =
\int_{\partial\Omega^\Gamma_\e}
|Du_\e|^{p} v\cdot\nu\, d\sigma 
\eeq
and
\beq
\label{e3}
\int_{\Omega^\Gamma_\e}f(u_\e)\, v\cdot Du_\e dx
=\int_{\Omega^\Gamma_\e}\sum_{i=1}^2v_i D_iF(u_\e)dx
=-\int_{\Omega^\Gamma_\e} F(u_\e)\div v\, dx.
\eeq
Finally, notice that
\begin{eqnarray}
\nonumber
\int_{\Omega^\Gamma_\e}\sum_{i,j=1}^2v_jD_{i,j}u_\e
|Du_\e|^{p-2}D_iu_\e dx
&=&
{1\over 2}\int_{\Omega^\Gamma_\e}\sum_{i,j=1}^2v_j|Du_\e|^{p-2}D_j|D_{i}u_\e|^2
dx
\\
\label{e4}
& = &
{1\over p}\int_{\Omega^\Gamma_\e}\sum_{j=1}^2v_j D_j|D u_\e|^p dx
\\
\nonumber
& = & 
{1\over p}\int_{\partial\Omega^\Gamma_\e}|Du_\e|^p v\cdot\nu d\sigma-
{1\over p}\int_{\Omega^\Gamma_\e}\div v |Du|^pdx.
\end{eqnarray}
Thus, (\ref{id}) follows combining (\ref{e1}), (\ref{e2}),
(\ref{e3}), (\ref{e4}).

\qed 

\begin{cor}
\label{C2.4}
Let $1<p<2$ and consider the vector field
$v\in\cC^1(\overline\Omega^\Gamma_{\bar\e_1},\R^2)$ defined by
(\ref{vf}).
Then, for all $\e\in(0,\bar\e_1)$, every solution $u_\e$ of the Dirichlet problem
\beq
\label{pP}
\div(|Du_\e|^{p-2}Du_\e)+f(u_\e)=0\qquad\mbox{ in
}\Omega^\Gamma_\e,\qquad u_\e=0\quad\mbox{ on
}\partial\Omega^\Gamma_\e
\eeq
satisfies the inequality
\beq
\label{le'}
0\le 
\left[ 1-{2\over p}+\left(1+{1\over p}\right)\mu(\e)\right]
\int_{\Omega^\Gamma_\e}|Du_\e|^pdx +
\int_{\Omega^\Gamma_\e}(\div v)\, F(u_\e)\, dx,
\eeq
where $\mu(\e)=\max\left\{
\left| {r[t\gamma_{\bar\e_1}''(t)\cdot
    N(t)]'\over 1-r\gamma_{\bar\e_1}''(t)\cdot
    N(t)}\right|\ :\ -\e\le r\le \e,\ t\in [a-\e,b+\e]
\right\}$.
\end{cor}

\no The proof follows directly from Lemma \ref{L2.1} and Lemma
\ref{L2.3}.

\vspace{2mm}

{\mbox {{\underline {\sf Proof of Theorem \ref{T2.0}}} \hspace{2mm}}}
First notice that assumption (\ref{f}) implies $f(0)=0$, so
the problem has the trivial solution $u\equiv 0$ for all $\e\in(0,\bar\e_1)$.
In order to prove that this solution is unique for $\e$ small enough,
for all $\e\in(0,\bar\e_1)$ let us consider a solution $u_\e$ of
problem (\ref{P}).
Taking into account condition (\ref{f}), from Lemma \ref{L2.1} and
Corollary \ref{C2.4} we obtain
\beq
\label{*}
0\le\left[1-{2\over p}+\left(1+{1\over p}\right)\mu(\e)\right]
  \int_{\Omega^\Gamma_\e}|Du_\e|^{p}dx+\left[2+\mu(\e)\right]
  {1\over q}\int_{\Omega^\Gamma_\e}u_\e f(u_\e)dx.
\eeq
Moreover, we have
\beq
\int_{\Omega^\Gamma_\e}u_\e f(u_\e)dx=\int_{\Omega^\Gamma_\e}|Du_\e|^{p}dx
\eeq
because $u_\e$ solves the Dirichlet problem (\ref{P}).
Therefore, (\ref{*})
implies
\beq
\label{le''}
0\le \left[1-{2\over p}+{2\over q}+\left(1+{1\over
      p}+{1\over
      q}\right)\mu(\e)\right]\int_{\Omega^\Gamma_\e}|Du_\e|^{p}dx.
\eeq
Since $\lim\limits_{\e\to 0}\mu(\e)=0$ (as follows from Lemma \ref{L2.2}) and
$1-{2\over p}+{2\over q}<0$ because $q>{2p\over 2-p}$, there exists
$\bar\e\in(0,\bar\e_1)$ such that 
\beq
1-{2\over p}+{2\over q}+\left(1+{1\over
      p}+{1\over
      q}\right)\mu(\e)<0
\qquad\forall \e\in(0,\bar\e).
\eeq
Therefore, if $\e\in(0,\bar\e)$ and $u_\e$ is a solution of the
Dirichlet problem in $\Omega^\Gamma_\e$, we must have
\beq
\int_{\Omega^\Gamma_\e}|Du_\e|^{p}dx=0
\eeq
that is $u_\e\equiv 0$ in $\Omega^\Gamma_\e$.
Thus, for all $\e\in(0,\bar\e)$, the problem has only the solution
$u\equiv 0$. 
\nopagebreak
\qed

\no Let us point out that Theorem \ref{T2.0} still holds if we replace
the smooth domain $\Omega=\Omega_\e^\Gamma$ by the piecewise smooth
domain
\beq
\cD_\e^\Gamma=\{\gamma(t)+r\, N(t)\ :\ t\in (a,b),\ |r|<\e\}.
\eeq
In this case, for the proof it is sufficient to apply the integral
identity given by Lemma \ref{L2.3} with the vector field $v$ defined
in (\ref{vf}) and to proceed as for the proof of Theorem \ref{T2.0},
taking into account that $v\cdot\nu\ge 0$ on $\partial\cD_\e^\Gamma$.


\sezione{The case $n>2$}


\no In Section 2, we proved that the Pohozaev type result for
star-shaped domains can be extended to a large class of contractible
non star-shaped domains of $\R^2$ so that the natural question arises
whether or not for $n=2$ this nonexistence result holds in all the
contractible domains.

\no Let us point out that the analogous question posed by Brezis for
$n\ge 3$ has negative answer, because there exist contractible domains
of $\R^n$ with $n\ge 3$, even arbitrarily close to non star-shaped
domains, such that the problem has nontrivial solutions.
This means that the existence of nontrivial solutions is related not
only to the topological but also to the metric properties of $\Omega$.

\no For example, for all $n\ge 2$, $\alpha\in\R$ and $d\in(0,1)$, let
us consider the bounded contractible domain $\cD_n^{\alpha,d}$ defined
by
\beq
\label{D}
\cD_n^{\alpha,d}=\left\{x=(x_1,\ldots,x_n)\in\R^n\ :\ 1-d<|x|<1+d,\
x_n<\alpha\left(\sum_{i=1}^{n-1}x_i^2\right)^{1/2}\right\}.
\eeq 
Then the following proposition holds (it is a particular case of some
existence and multiplicity results obtained in
\cite{MPaihp04,MPcvpde06,MPaihp06,Pmm89,P4092,P94,Ptmna96,Pd98})

\begin{prop}
\label{P3.1}
Let $n\ge 3$, $d\in(0,1)$ and $q\ge {2n\over n-2}$.

\no Then, there exists $\bar\alpha\in\R$ such that, if
$\alpha\ge\bar\alpha$, the Dirichlet problem 
\beq
\label{star}
\Delta u+|u|^{q-2}u=0\quad\mbox{ in }\cD_n^{\alpha,d},\qquad
u=0\quad\mbox{ on }\partial\cD_n^{\alpha,d}
\eeq
has positive and sign changing solutions.
Moreover, as $\alpha\to +\infty$, these solutions tend to 0 and their
number tends to infinity.
\end{prop}

\no On the contrary, if $n=2$, $1<p<2$ and condition (\ref{f}) holds
for $q>{2p\over 2-p}$, the result obtained in Section 2 guarantees in
particular that there exists $\bar\e\in(0,1)$ such that Problem
(\ref{P}) with $\Omega=\cD_2^{\alpha,d}$ has only the trivial solution
$u\equiv 0$ for all the pairs $(\alpha,d)$ such that $\alpha\in\R$ and
$d\in(0,\bar\e)$. 

\no Notice that the contractible domain $\cD_n^{\alpha,d}$ tends as
$\alpha\to+\infty$ to the annulus $A_n^d=\{x\in\R^n\ :\
1-d<|x|<1+d\}$, which is non contractible in itself.
Thus, in the contractible domains $\cD_n^{\alpha,d}$, which are
arbitrarily close to non contractible domains for $\alpha$ large enough,
there exists only the trivial solution $u\equiv 0$ if $n=2$, while
there exist many nontrivial solutions if $n>2$.

\no Moreover, we have the following proposition where we gather some
existence and multiplicity results that are particular cases of more
general results obtained in \cite{MPcras02,MPaihp04,MPcvpde06} etc..

\begin{prop}
\label{P3.2}
Let $n\ge 3$, $\alpha>0$, $d\in(0,1)$ and consider the domain
$\cD_n^{\alpha,d}$ defined in (\ref{D}).

\no Then there exists $\bar q\ge {2n\over n-2}$ and $\bar\e>0$ such
that problem (\ref{P}) with $\Omega=\cD_n^{\alpha,d}$ has solutions
for all $q\ge\bar q$ and for all $q\in\left({2n\over n-2},{2n\over
    n-2}+\bar\e\right)$.

\no Moreover, these solutions tend to 0 as $q\to\infty$ and $q\to  
{2n\over n-2}$, while their number tends to infinity.
\end{prop}

\no Notice that the domain $\cD_n^{\alpha,d}$ is non star-shaped if
$\alpha>0$ while if $\alpha<0$ it is star-shaped for $d$ close to 1,
so the problem has only the trivial solution $u\equiv 0$ (this means
that the result given in Proposition \ref{P3.2} is sharp for what
concerns the assumption on $\alpha$).

\no Next theorem (which extends Theorem \ref{T2.0} to the case $n\ge
3$) shows that, as for $n=2$, also for $n\ge 3$ there exist suitable
contractible non star-shaped domains, even arbitrarily close to non
contractible domains, such that the problem has only the trivial
solution $u\equiv 0$ (see also \cite{MPtub} for related results).
Taking into account Proposition \ref{P3.1}, it is clear that these
domains and the contractible domains $\cD_n^{\alpha,d}$ with $\alpha$
large must have quite different geometrical properties (as we explain
in Remark \ref{RB}).

\begin{teo}
\label{T3.0}
Let $n>2$, $1<p<n$, and assume that condition (\ref{f}) holds for a
suitable $q>{np\over n-p}$.
Let $\Gamma$ and $\Omega^\Gamma_\e$ be as in Theorem \ref{T2.0} and
consider the domains $\Omega=\widetilde\Omega^\Gamma_{\e,s}$ defined in
(\ref{Otilde}).
Then, there exists $\tilde\e>0$ such that the Dirichlet problem
(\ref{P}) has only the trivial solution $u\equiv 0$ for all the pairs $(\e,s)$
such that $\e\in(0,\tilde\e)$ and $s>0$.
\end{teo}

\no In order to prove Theorem \ref{T3.0}, we proceed as in Section 2,
but now we use the vector field $\tilde v=(\tilde
v_1,\ldots, \tilde v_n)\in\cC^1(\widetilde\Omega^\Gamma_{\e,s},\R^n)$
defined by 
$$
\tilde v(\gamma_{\bar\e_1}(t)+rN(t),y_1,\ldots,y_{n-2})=
(tT(t)[1-r\,
\gamma_{\bar\e_1}''(t)\cdot N(t)]+rN(t), y_1,\ldots,y_{n-2})
$$
\beq
\phantom{***********}\forall t\in[a-{\bar\e_1},
b+{\bar\e_1}],\ \forall r\in[-{\bar\e_1},{\bar\e_1} ],\ \forall
(y_1,\ldots,y_{n-2})\in\R^{n-2}. 
\eeq
Then, Lemma \ref{L2.1} has to be modified as follows.

\begin{lemma}
\label{L3.1}
If $\gamma$ is as in Lemma \ref{L2.1}, the vector field $\tilde v$
satisfies
\begin{itemize}
\item[{\em a)}]$\tilde v\cdot\tilde \nu>0$ on $\partial\widetilde\Omega^\Gamma_{\e,s}$
  $\forall\e\in(0,\bar\e_1)$ $\forall s>0$ where $\tilde\nu$ denotes
  the outward normal to $\partial \widetilde\Omega^\Gamma_{\e,s}$; 
\item[{\em b)}] $\div
  \tilde v(\gamma_{\bar\e_1}(t)+rN(t),y)=n-{r[t\gamma_{\bar\e_1}''(t)\cdot
    N(t)]'\over 1-r\gamma_{\bar\e_1}''(t)\cdot
    N(t)}$ $\forall t\in[a-{\bar\e_1},
b+{\bar\e_1}]$, $\forall r\in[-{\bar\e_1},{\bar\e_1} ]$,
$\forall y\in\R^{n-2}$;
\item[{\em c)}]
$d\tilde v(\gamma_{\bar\e_1}(t)+rN(t),y)[\xi,\psi]\cdot(\xi,\psi)=
\left[1-{r[t\gamma_{\bar\e_1}''(t)\cdot
    N(t)]'\over 1-r\gamma_{\bar\e_1}''(t)\cdot
    N(t)}\right]\,\xi^2_T(t)+\xi^2_N(t)+|\psi|^2$ $\forall
t\in[a-{\bar\e_1}, 
b+{\bar\e_1}]$, $\forall r\in[-{\bar\e_1},{\bar\e_1} ]$, $\forall
y\in\R^{n-2}$, $\forall\xi\in\R^2$, $\forall\psi\in\R^{n-2}$, where 
$\xi_T(t)=\xi\cdot T(t)$
and 
$\xi_N(t)=\xi\cdot N(t)$.
\end{itemize}
\end{lemma}

\proof Property $(a)$ follows directly from the definition of $\tilde
v$ and $ \widetilde\Omega^\Gamma_{\e,s}$ (which is a piecewise smooth
domain).

\no The proof of $(b)$ and $(c)$ is as in Lemma \ref{L2.1} taking also
into account that 
\beq
d\tilde v(\gamma_{\bar\e_1}(t)+rN(t),y)[0,\psi]=(0,\psi)
\eeq
as one can verify by direct computation.

\qed

\no Lemma \ref{L2.2}, Lemma \ref{L2.3}, Corollary \ref{C2.4} and their proofs
require only obvious modifications that take into account Lemma
\ref{L3.1}.
In particular, the inequality (\ref{le'}) in Corollary \ref{C2.4}
becomes now
\beq
0\le 
\left[ 1-{n\over p}+\left(1+{1\over p}\right)\mu(\e)\right]
\int_{\widetilde\Omega^\Gamma_{\e,s}}|D\tilde u_\e|^pdx +
\int_{\widetilde\Omega^\Gamma_{\e,s}}(\div \tilde v) \, F(\tilde u_\e)\,dx
\eeq
for all solutions $\tilde u_\e$ of the Dirichlet problem in the domain 
$ \widetilde\Omega^\Gamma_{\e,s}$ and, as a consequence, the
inequality (\ref{le''}) becomes
\beq
0\le \left[1-{n\over p}+{n\over q}+\left(1+{1\over
      p}+{1\over
      q}\right)\mu(\e)\right]\int_{\widetilde\Omega^\Gamma_{\e,s}}|D\tilde
u_\e|^{p}dx.
\eeq
Then, since $1-{n\over p}+{n\over q}<0$ for $q>{np\over n-p}$, the
proof of Theorem \ref{T3.0} may be carried out following the same
procedure as in the proof of Theorem \ref{T2.0}. 

\begin{rem}{\em
\label{RA}
Notice that in the domains $\widetilde\Omega_{\e,s}^\Gamma$ arising in
Theorem \ref{T3.0} only $\e$ is required to be small while $s$ may be
arbitrarily large.
This means that these domains are thin only in one dimension (while
the domains considered in \cite{MPtub} are thin in $n-1$ dimensions).

\no Moreover, taking into account the definition of the vector field
$\tilde v$ used in the proof of Theorem \ref{T3.0}, one can verify by
direct computation that this theorem still holds if the domains
$\widetilde\Omega_{\e,s}^\Gamma$ are replaced by the more general
domains
\beq
\label{sigma}
\widetilde\Omega_{\e,s}^{\Gamma,\Sigma}=\{\gamma(t)+r\, N(t)+(0,0,y)\ :\
t\in(a,b),\ |r|<\e,\ y\in\R^{n-2},\ |y|<s,\ (r,y)\in\Sigma\}
\eeq
where $\Sigma$ is a domain of $\R^{n-1}$, star-shaped with respect to
the origin.
In particular, Theorem \ref{T3.0} holds for the domains
$\Omega=\widetilde\cD_{\e,s}^\Gamma$ defined by
\beq
\widetilde\cD_{\e,s}^\Gamma=\{\gamma(t)+r\, N(t)+(0,0,y)\ :\
t\in(a,b),\ |r|<\e,\ y\in\R^{n-2},\ |y|<s\}
\eeq
(that are obtained, for example, when $\Sigma=\R^{n-1}$ in (\ref{sigma})).

\no In fact, for the proof, we need only to verify that property $(a)$
in Lemma \ref{L3.1} still holds if $\widetilde\Omega_{\e,s}^\Gamma$ is
replaced by $\widetilde\Omega_{\e,s}^{\Gamma,\Sigma}$.

\no Notice that the class of the domains  $\widetilde\Omega_{\e,s}^{\Gamma,\Sigma}$
includes also domains of the form 
\beq
\widetilde\cD_{d,s}^\alpha \{(x_1,x_2,y)\in\R^n\ :\
(x_1,x_2)\in\cD_2^{\alpha,d},\ y\in\R^{n-2},\ |y|<s\}
\eeq
that are obtained when in (\ref{sigma}) $\gamma$ is an arc of
circumference and $\Sigma=\R^{n-1}$.

}\end{rem}

\begin{rem}{\em
\label{RB}
In order to explain the sense of these results in the framework of the
study of the effect of the domain shape on existence and nonexistence
of nontrivial solutions for nonlinear elliptic problems with critical
and supercritical growth, let us recall that the results obtained in
\cite{MPaihp04,MPcvpde06,MPaihp06,Pmm89,P4092,P94,Ptmna96,Pd98}
suggest that the number of nontrivial solutions for these problems is
related to the property that the domain $\Omega$ can be obtained by
removing a subset of small capacity from a domain having different
$k$-dimensional homology group with $k\ge 2$.

\no Thus, the existence and multiplicity result in the domains of the
form $\cD_n^{\alpha,d}$ with $n\ge 3$ and $\alpha$ large enough, given
by Proposition \ref{P3.1}, is related to the fact that the
contractible domain $\cD_n^{\alpha,d}$ tends as $\alpha\to +\infty$ to
the annulus  $A_n^d$ which has different $(n-1)$-dimensional homology
group (with $n-1\ge 2$) and the capacity of
$A_n^d\setminus\cD_n^{\alpha,d}$ tends to 0 as $\alpha\to +\infty$.

\no On the contrary, the contractible domains $\cD_2^{\alpha,d}$ and
$\widetilde\cD_{d,s}^\alpha$ (see Remark \ref{RA}) tend as $\alpha\to
+\infty$ to non contractible domains where only the 1-dimensional
homology group is nontrivial; moreover, these domains do not differ
from their limit domains by sets whose capacities tend to 0 as
$\alpha\to +\infty$.
These facts explain the deep reason of the nonexistence results given
by Theorems \ref{T2.0} and \ref{T3.0}.

}\end{rem}

\vspace{5mm}


{\small {\bf Acknowledgement}. The authors have been supported by the ``Gruppo
Nazionale per l'Analisi Matematica, la Probabilit\`a e le loro
Applicazioni (GNAMPA)'' of the {\em Istituto Nazio\-nale di Alta Matematica
(INdAM)} - Project: Equazioni di Schrodinger nonlineari: soluzioni con
indice di Morse alto o infinito. 

The second author acknowledges also the MIUR Excellence Department
Project awarded to the Department of Mathematics, University of Rome
Tor Vergata, CUP E83C18000100006. 
}



\end{document}